\documentclass{amsart}
\usepackage{amsthm}
\usepackage{paralist}
\usepackage{pictexwd,dcpic}
\usepackage{url}
\usepackage{texdraw}
\usepackage{graphics}
\usepackage{graphicx}
\usepackage[english]{babel}
\usepackage{amssymb,amsmath,latexsym} 
\usepackage{epic,eepic}
\usepackage[latin1]{inputenc}
\numberwithin{equation}{section}
\def\neu#1{{\bf #1}}
\def\KK{\mathbb{K}}
\def\NN{\mathbb{N}}

\def\bb#1{\mathbb{#1}}
\def\frak#1{\mathfrak{#1}}

\def\cal#1{\mathcal{#1}}
\def\tru{\triangle}
\def\trd{\bigtriangledown}
\def\tu{^\tru}
\def\td{^\trd}
\def\wdl{weakly dicomplemented lattice}
\def\I{\mathop{\mbox{\rm I}}}

\begin{document}

\theoremstyle{plain}
\newtheorem{definition}{Definition}[section]
\newtheorem{theorem}{Theorem}[section]
\newtheorem{proposition}[theorem]{Proposition}
\newtheorem{lemma}[theorem]{Lemma}
\newtheorem{corollary}[theorem]{Corollary}
\newtheorem{example}{Example}[section]
\newtheorem{remark}{Remark}[section]

\title{On the isomorphism problem of concept algebras}
\author{L\'eonard Kwuida}
\address[L. Kwuida]{Zurich University of Applied Sciences\\ School of Engineering\\ Technikumstrasse 9 \\ CH-8401 Winterthur\\ Switzerland}
\email{kwuida@gmail.com}

\author{Hajime Machida}
\address[H. Machida]{Hitotsubashi University\\ Department of Mathematics\\  2-1 Naka, Kunitachi, Tokyo 186-8601, Japan}
\email{machida@math.hit-u.ac.jp}
\thanks{This paper is an extended version of \cite{KM08}}
\begin{abstract}
Weakly dicomplemented lattices are bounded lattices equipped with two unary operations to encode a negation on {\it concepts}. They have been introduced to capture the equational theory of concept algebras~\cite{Wi00}. They generalize Boolean algebras. Concept algebras are concept lattices, thus complete lattices, with a weak negation and a weak opposition. A special case of the representation problem for weakly dicomplemented lattices, posed in \cite{Kw04}, is whether complete {\wdl}s are isomorphic to concept algebras. In this contribution we give a negative answer to this question (Theorem~\ref{T:main}). We also provide a new proof of a well known result due to M.H. Stone~\cite{St36}, saying that {\em each Boolean algebra is a field of sets} (Corollary~\ref{C:Stone}). Before these, we prove that the boundedness condition on the initial definition of {\wdl}s (Definition~\ref{D:wdl}) is superfluous (Theorem~\ref{T:wcl}, see also \cite{Kw09}). 
\end{abstract}

\maketitle
\section{Introduction}
Formal concept analysis (FCA) started in the 80ies from an attempt to restructure lattice theory by Rudolf Wille~\cite{Wi82}. FCA is based on the formalization of the notions of ``concept'' and ``concept hierarchy''. In traditional philosophy a concept is defined by its extent and its intent: the extent contains all entities belonging to the concept, and the intent contains all properties satisfied by exactly all entities of the concept. The concept hierarchy states that ``a concept is more general if it contains more objects, or equivalently, if its intent is smaller''. The set of all concepts of a ``context'' with its concept hierarchy forms a complete poset called {\em concept lattice}. Based on ordered structures, FCA provides a nice formalism for knowledge management and retrieval. It has developed rapidly and now stands as a research area on its own, and has been applied in many fields. For displaying knowledge FCA offers several techniques, among them the line diagrams (visualization) and  the implication theory (logical description of the information~\cite{GD86,GW99a}). 

In his project to extend FCA to a broader field called {\em Contextual Logic}, Rudolf Wille needed to formalize a conceptual negation. The problem of negation is surely one of the oldest problems of the scientific and philosophic community, and still attracts the attention of many researchers (see \cite{Hl89,Wa96}). Several types of logic have been introduced, according to the behavior of the corresponding negation. To develop a contextual logic, one of the starting points is that of Boolean algebras, which arise from the encoding of the operations of human thought by George Boole~\cite{Bo54}. Is there a natural generalization of Boolean algebras to concept lattices?  Boolean Concept Logic aims to develop a mathematical theory for Logic, based on concept as unit of thought, as a generalization of that developed by George Boole in~\cite{Bo54}, based on signs and classes. The main operations of human mind that Boole encoded are conjunction, disjunction, universe, "nothing" and ``negation''. 

The set of all formal concepts of a given formal context forms a complete
lattice. Therefore, apart from the negation, the operations encoded
by Boole are without problem encoded by lattice operations. To encode a
negation Wille followed Boole's idea, and suggested many candidates, among them 
a weak negation (taking the concept generated by the complement of the extent) and a weak opposition (taking the concept generated by the complement of the intent)~\cite{Wi00}. This approach is driven by the wish to have a negation as an internal operation on concepts\footnote{Other approaches have to relax the definition of concept. These are preconcepts, semiconcepts and protoconcepts. They have been investigated by Rudolf Wille and coworkers for example in \cite{Wi00,HLSW,VW05,BW06},\dots. In \cite{DE98}, there is another proposition to get negation on lattices.}. The concept lattice together with these operations is called concept algebra. Expressing a negation in information science and knowledge systems
can be very helpful, in particular while dealing with incomplete information (see for example \cite{MNR08,Pr06,BH05,Fe06}). In the absence of a Boolean negation, weak negation and weak opposition would offer an alternative. In this case concept algebras and  weakly dicomplemented lattices (see below) would replace powerset algebras and Boolean algebras respectively.

For abstracting concept algebras, weakly dicomplemented lattices have been introduced. Those are lattices with two unary operations that satisfy some equations known to hold in concept algebras. The main problem we address in this paper is when a weakly dicomplemented lattice is isomorphic to a concept algebra. Characterizing concept algebras remains an open problem, but substantial results are obtained, especially in the finite case~\cite{Kw04,GK07}. The rest of this contribution is divided as follows: in Section~\ref{S:wdicomp} we introduce some formal definitions, give a characterization of {\wdl}s and present several constructions of {\wdl}s. Section~\ref{S:wdl+neg} shows why weakly dicompelemented lattices are considered as a generalization of Boolean algebras. In Section~\ref{S:main} we prove that completeness is not enough to get {\wdl}s isomorphic to concept algebras. We end with a new proof of the representation of Boolean algebras by fields of sets. 

\section{Weak dicomplementation.}\label{S:wdicomp}
\begin{definition}\label{D:wdl}
A \neu{weakly dicomplemented lattice} is a bounded lattice $L$ equipped with two unary operations $\tu$ and $\td$ called \neu{weak complementation} and \neu{dual weak complementation}, and satisfying for all $x,y\in L$ the following equations\footnote{Note that $x^{\tru\tru}\le x\iff x^{\tru\tru}\vee x=x$ and $x^{\trd\trd}\ge x\iff x^{\trd\trd}\wedge x =x$; thus conditions (1) and (1') can be written as equations. For (2) and (2') we have $x\le y\implies x\tu\ge y\tu$
is equivalent to $(x\wedge y)\tu\wedge y\tu= y\tu$ and $x\le y\implies x\td\ge y\td$ is equivalent to $(x\wedge y)\td\wedge y\td =y\td$.}: 
\medskip\par\noindent
\begin{minipage}{.45\textwidth}
\begin{enumerate}
\item[(1)] $x^{\tru\tru}\le x$,
\item[(2)] $x\le y\implies x\tu\ge y\tu$,
\item[(3)] $(x\wedge y)\vee(x\wedge y\tu)=x$,
\end{enumerate}
\end{minipage}\hfill
\begin{minipage}{.45\textwidth}
\begin{enumerate}
\item[(1')] $x^{\trd\trd}\ge x$,
\item[(2')] $x\le y\implies x\td\ge y\td$,
\item[(3')] $(x\vee y)\wedge(x\vee y\td)=x$.
\end{enumerate}
\end{minipage}\medskip\par\noindent
\noindent We call $x\tu$ the \neu{weak complement} of $x$ and $x\td$ the \neu{dual weak complement} of $x$. The pair $(x\tu,x\td)$ is called the \neu{weak dicomplement} of $x$ and the pair $(\tu,\td)$ a  \neu{weak dicomplementation} on $L$. The structure $(L,\wedge,\vee,\tu,0,1)$ is  called a \neu{weakly complemented lattice} and $(L,\wedge,\vee,\td,0,1)$ a \neu{dual weakly complemented lattice}.
\end{definition}
The following properties are easy to verify: \qquad
(i)\quad $x\vee x\tu=1$, \quad(ii)\quad {$x\wedge x\td=0$}, \newline
(iii)\quad $0\tu=1=0\td$,\quad (iv)\quad $1\tu=0=1\td$,\quad (v)\quad $x\td\leq x\tu$,\quad (vi)\quad $(x\wedge y)\tu=x\tu\vee y\tu$, \newline
 (vii)\quad $(x\vee y)\td=x\td\wedge y\td$,\quad (viii)\quad $x\sp{\tru\tru\tru}=x\tu$,\quad (ix)\quad $x\sp{\trd\trd\trd}=x\td$ and  \newline 
(x)\quad $x{\tu}{\td}\leq x{\tu}{\tu}\leq x \leq x{\td}{\td}\leq x{\td}{\tu}$. 
\begin{example}\label{E:ba}
\hfill
\begin{enumerate}
\item[(a)] The natural examples of {\wdl}s are Boolean algebras. For a Boolean algebra $(B,\wedge,\vee,\bar{\phantom{a}},0,1)$, the algebra $(B,\wedge,\vee,\bar{\phantom{a}},\bar{\phantom{a}},0,1)$ (complementation duplicated, i.e. $x\tu := \bar{x}=: x\td$) is a weakly dicomplemented lattice.
\item[(b)]  Each bounded lattice can be endowed with a {\bf trivial weak dicomplementation} by de\-fi\-ning $(1,1),\,(0,0)$ and $(1,0)$ as the dicomplement of $0,\ 1$ and of each ${x\not\in\{0,1\}}$, respectively.
\end{enumerate}
\end{example}
\begin{theorem}\label{T:wcl}
Weakly complemented lattice are exactly the nonempty lattices with an additional unary operation $\tu$ that satisfy the equations (1)--(3) in Definition~\ref{D:wdl}. 
\end{theorem}
Of course, weakly complemented lattices satisfy the equations (1)--(3) in Definition~\ref{D:wdl}. So what we should prove is that all non empty lattices satisfying the equations $(1)-(3)$ are bounded. 
\begin{proof}
Let $L$ be a nonempty lattice satisfying the equations (1)--(3). For an element $x\in L$, we set ${1:=x\vee x\tu}$ and $0:=1\tu$. We are going to prove that $1$ and $0$ are respectively the greatest and lowest element of $L$. Let $y$ be an arbitrary element of $L$. We have 
$$1\geq y\wedge 1=y\wedge(x\vee x\tu)\geq(y\wedge x)\vee(y\wedge x\tu)=y, \quad \text{ by (3)}.$$
 Thus $x\vee x\tu$ is the greatest element of $L$. Of course, if $L$ was equipped with a unary operation $\td$ satisfying the equation (1')--(3') we could use the same argument as above to say that $x\wedge x\td$ is the smallest element of $L$. Unfortunately we have to check that $0:=1\tu$ is less than every other element of $L$. So let $y\in L$. We want to prove that $0\leq y$. Note that 
\[
(y\wedge y\tu)\tu\geq y\tu\vee y^{\tru\tru}=1.\]
Thus $(y\wedge y\tu)\tu=1$.
 For an arbitrary element $z$ of $L$, we have 
\[0\wedge z=1\tu\wedge z=(y\wedge y\tu)^{\tru\tru}\wedge z\leq y\wedge y\tu\wedge z\leq y\wedge z\]
and
\[0\wedge z\tu=1\tu\wedge z\tu=(y\wedge y\tu)^{\tru\tru}\wedge z\tu\leq y\wedge y\tu\wedge z\tu\leq y\wedge z\tu.\]
Henceforth $0=(0\wedge z)\vee(0\wedge z\tu)\leq (y\wedge z)\vee(y\wedge z\tu)=y$.  
\end{proof}  
\begin{remark}
In Universal Algebra~\cite{BS81}, one should care about the signature while defi\-ning an algebra. By Theorem~\ref{T:wcl} we can choose between $(\wedge,\vee,\tu)$ and $(\wedge,\vee,\tu,0,1)$ as signature for weakly complemented lattices. Let $\mathbb{V}_1$ be the variety of algebras $(L,\wedge,\vee,\tu)$ of type $(2,2,1)$ such that $(L,\wedge,\vee)$ is a lattice satisfying the equations (1)--(3) in Definition~\ref{D:wdl}, and $\mathbb{V}_2$ the variety of algebras of type $(2,2,1,0,0)$ such that $(L,\wedge,\vee,0,1)$ is a bounded lattice satisfying the equations (1)--(3) in Definition~\ref{D:wdl}. Then an algebra with the empty set as carrier set belongs to $\mathbb{V}_1$, but not to $\mathbb{V}_2$.
Any non empty substructure of an algebra of $\mathbb{V}_1$ is a substructure of the corresponding algebra in $\mathbb{V}_2$ and vice versa. 
 Any map that is a morphism between nonempty algebras of $\mathbb{V}_1$ is also a morphism between algebras of $\mathbb{V}_2$ and vice-versa. Hence, there is no big diffe\-ren\-ce is considering one signature instead of another. Here we will keep the signature $(\wedge,\vee,\tu,0,1)$, to indicate ``contradiction'' and ``universe''. 
\end{remark}
\begin{definition}
Let $(P,\leq)$ be a poset and $f:P\to P$ be a map. $f$ is a \neu{closure operator} on $P$ if 
\[
x\leq f(y)\iff f(x)\leq f(y), \quad \text{ for all } x,y\in P.
\] 
This is equivalent to 
\[x\leq f(x),\quad x\leq y\implies f(x)\leq f(y)\quad \text{ and }\quad f(f(x))=f(x), \quad \text{ for all } x,y\in P.\]
 Usually we will write a closure operator on a set $X$ to mean a closure operator on the po\-wer\-set $(\cal{P}(X),\subseteq)$ of $X$.  
Dually, $f$ is a \neu{kernel operator} on $P$ if 
\[
x\geq f(y)\iff f(x)\geq f(y), \quad \text{ for all } x,y\in P.
\]
 As above, we say that $f$ is a kernel operator on $X$ to mean a kernel operator on $(\cal{P}(X),\subseteq)$.
\end{definition}
For a \wdl\ $(L,\wedge,\vee,\tu,\td,0,1)$, the maps $x\mapsto x\sp{\tru\tru}$ and $x\mapsto x\sp{\trd\trd}$ are resp. kernel and closure operators on $L$. If $f$ is a closure operator (resp. a kernel operator) on a lattice $L$, then $f(L)$ (with the induced order) is a lattice. Recall that for any closure operator $h$ on $L$ we have 
\[
h(h(x)\wedge h(y))=h(x)\wedge h(y)\quad \text{ as well as } \quad h(h(x)\vee h(y))=h(x\vee y);
\]
 dually, for any kernel operator $k$ on $L$ we have 
\[
k(k(x)\wedge k(y))=k(x\wedge y)\quad \text{ and }\quad k(k(x)\vee k(y))=k(x)\vee k(y).
\]
 We denote by $P^d$ the dual poset of $(P,\leq)$, i.e. $P^d:=(P,\geq)$. Then $f$ is a kernel operator on $P$ iff $f$ is a closure operator on $P^d$. 
\begin{proposition}\label{L:hk}
Let $h$ be a closure operator on a set $X$ and $k$ a kernel operator on a set $Y$. For $A\subseteq X$ and $B\subseteq Y$ define $A\sp{\tru_h}:=h(X\setminus A)$ and $B\sp{\trd_k}:=k(Y\setminus B)$.
\begin{itemize}
\item[(i)] $({h}\cal{P}(X),\cap,\vee^h,\sp{\tru_h},{h}\emptyset,X)$, with $A_1\vee^h A_2:=h(A_1\cup A_2)$, is a weakly complemented lattice.
\item[(i')]  $({k}\cal{P}(Y),\wedge_k,\cup,\sp{\trd_k},\emptyset,{k}Y)$, with $B_1\wedge_k B_2:=k(B_1\cap B_2)$, is a dual weakly complemented lattice.
\item[(ii)] If ${h}\cal{P}(X)$ is isomorphic to ${k}\cal{P}(Y)$, then $h$ and $k$ induce weakly dicomplemented lattice structures on ${h}\cal{P}(X)$ and on ${k}\cal{P}(Y)$ that are extensions of those in $(i)$ and $(i')$ above respectively.
\end{itemize}
\end{proposition}
\begin{proof}
For (i), let $h$ be a closure operator on X; $({h}\cal{P}(X),\cap,\vee^h,{h}\emptyset,X)$ is a bounded lattice. So we should only check that the equations $(1)-(3)$ in Definition~\ref{D:wdl} hold.  For $x\in{h}\cal{P}(X)$, we have $x\sp{\tru\tru}={h}(X\setminus{h}(X\setminus x))\subseteq{h}(X\setminus(X\setminus x))=h(x)=x$, and $(1)$ is proved. For $x_1\leq x_2$ in ${h}\cal{P}(X)$, we have $x_1\subseteq x_2$ and ${h}(X\setminus x_1)\supseteq{h}(X\setminus x_2)$, and $(2)$ is proved. Now we consider $x,y\in {h}\cal{P}(X)$. Trivially $(x\cap y)\vee^h(x\cap y\sp{\tru_h})\leq x$. In addition, 
\[(x\cap y)\vee^h(x\cap y\sp{\tru_h})=(x\cap y)\vee^h(x\cap h(X\setminus y))=h((x\cap y)\cup(x\cap h(X\setminus y)))\]
\[\supseteq h((x\cap y)\cup(x\cap(X\setminus y)))= h(x)= x.\]
$(i')$ is proved similarly. 

For $(ii)$ we will extend the structures of $(i)$ and $(i')$ to get {\wdl}s. By $(i)$, $({h}\cal{P}(X),\cap,\vee^h,\sp{\tru_h},{h}\emptyset,X)$ is a weakly complemented lattice. Let $\varphi$ be an isomorphism from ${h}\cal{P}X$ to ${k}\cal{P}Y$. We define $\sp{\trd_\varphi}$ on ${h}\cal{P}(X)$ by: 
$$x\sp{\trd_\varphi}:=\varphi^{-1}(\varphi(x)\sp{\trd_k}).$$
 Then 
$$x\sp{\trd_\varphi\trd_\varphi}=\left(\varphi^{-1}\left(\varphi(x)\sp{\trd_k}\right)\right)\sp{\trd_\varphi} = \varphi^{-1}\left(\varphi\left(\varphi^{-1}\left(\varphi(x)\sp{\trd_k}\right)\right)\sp{\trd_k} \right)=\varphi^{-1}\left(\varphi(x)\sp{\trd_k\trd_k} \right),$$
$\text{and }x\sp{\trd_\varphi\trd_\varphi}\geq\varphi^{-1}(\varphi(x))=x$. 
 For $x\leq y$ in ${h}\cal{P}X$ we have 
$\varphi(x)\leq\varphi(y)$ implying 
\[
\varphi(x)\sp{\trd_k}\geq\varphi(y)\sp{\trd_k}\quad \text{ and }\quad x\sp{\trd_\varphi}=\varphi^{-1}(\varphi(x)\sp{\trd_k})\geq\varphi^{-1}(\varphi(y)\sp{\trd_k})=y\sp{\trd_\varphi}.
\]
 For $x,y$ in ${h}\cal{P}X$, we have
\begin{eqnarray*}
(x\vee y)\wedge(x\vee y\sp{\trd_\varphi})&=&(x\vee y)\wedge(x\vee \varphi^{-1}(\varphi(y)\sp{\trd_k}))\\
&=&\varphi^{-1}\left((\varphi(x)\vee \varphi(y))\wedge(\varphi(x)\vee \varphi(y)\sp{\trd_k})\right)\\
&=&\varphi^{-1}(\varphi(x))=x.
\end{eqnarray*}
 Therefore 
$({h}\cal{P}(X),\cap,\vee^h,\sp{\tru_h},\sp{\trd_\varphi},{h}\emptyset,X)$ is a \wdl. Similarly \linebreak $({k}\cal{P}(Y),\wedge^k,\cup,\sp{\tru_\varphi},\sp{\trd_k},\emptyset,{k}Y)$ with $x\sp{\tru_\varphi}:=\varphi(\varphi^{-1}(x)\sp{\tru_h})$ is a \wdl.  
\end{proof}
\begin{proposition}\label{L:combinehk}
Let $h$ be a closure operator on $X$ and $k$ a kernel operator on $Y$ such that ${h}\cal{P}(X)$ is isomorphic to ${k}\cal{P}(Y)$. Let $\varphi$
be an isomorphism from ${h}\cal{P}(X)$ to ${k}\cal{P}(Y)$. 
We set 
\[L:=\{(x,y)\in {h}\cal{P}(X)\times {k}\cal{P}(Y)\mid y=\varphi(x)\}.
\]
 $L$ has a \wdl\ structure induced by $h$ and $k$.
\end{proposition}
\begin{proof}
By Lemma~\ref{L:hk} $({h}\cal{P}(X),\cap,\vee^h,\sp{\tru_h},{h}\emptyset,X)$ is a weakly complemented lattice and  \linebreak $({k}\cal{P}(Y),\wedge^k,\cup,\sp{\trd_k},\emptyset,{k}Y)$ a dual weakly complemented lattice. 
 For every $y\in {k}\cal{P}(Y)$ there is a unique $x\in{h}\cal{P}(X)$ such that $y=\varphi(x)$. For $(a,b)$ and $(c,d)$ in $L$, we have $a\leq c\iff b\leq d$. 
 We define a relation $\leq$ on $L$ by: $a\leq c\iff:(a,b)\leq(c,d):\iff b\leq d$. Then 
\[
{h}\cal{P}(X)\stackrel{\pi_1}{\cong}L\stackrel{\pi_2}{\cong}{k}\cal{P}(Y)
\]
 where $\pi_i$ is the $i^{th}$ projection. Thus $(L,\leq)$ is a bounded lattice. Moreover 
\[
(a,b)\wedge(c,d)=(a\cap c, \varphi(a\cap c))\quad \text{ and } \quad (a,b)\vee(c,d)=(\varphi^{-1}(b\cup d),b\cup d).
\]
 For $(x,y)\in L$, we define $(x,y)\tu:=(x\sp{\tru_h},\varphi(x\sp{\tru_h}))$ and $(x,y)\td:=(\varphi^{-1}(y\sp{\trd_k}),y\sp{\trd_k})$. We claim that $(L,\wedge,\vee,\tu,\td,0,1)$ is a \wdl. In fact,
\[(x,y)\sp{\tru\tru}=(x\sp{\tru_h},\varphi(x\sp{\tru_h}))\tu=(x\sp{\tru_h\tru_h},\varphi(x\sp{\tru_h\tru_h}))\leq(x,\varphi(x))=(x,y).\]
If $(x,y)\leq (z,t)$ in $L$, we have $x\leq z$ and $y\leq t$, implying $x\sp{\tru_h}\geq z\sp{\tru_h}$ and \newline  ${\varphi(x\sp{\tru_h})\geq\varphi(z\sp{\tru_h})}$; thus $(x,y)\tu=(x\sp{\tru_h},\varphi(x\sp{\tru_h}))\geq(z\sp{\tru_h},\varphi(z\sp{\tru_h}))=(z,t)\tu$. These prove (1) and (2) of Definition~\ref{D:wdl}. It remains to prove (3). Let $(x,y)$ and $(z,t)$ in $L$;
{ 
\begin{eqnarray*}
((x,y)\wedge(z,t))\vee((x,y)\wedge(z,t)\tu) &=& (x\cap z, \varphi(x\cap z))\vee((x,y)\wedge(z\sp{\tru_h},\varphi(z\sp{\tru_h})))\\
&=&(x\cap z, \varphi(x\cap z))\vee(x\cap z\sp{\tru_h},\varphi(x\cap z\sp{\tru_h}))\\
&=&(\varphi^{-1}(\varphi(x\cap z)\cup\varphi(x\cap z\sp{\tru_h})), \varphi(x\cap z)\cup\varphi(x\cap z\sp{\tru_h}))\\
&=& ((x\cap z)\vee^h(x\cap z\sp{\tru_h}),\varphi((x\cap z)\vee^h(x\cap z\sp{\tru_h})))\\
&=& (x,\varphi(x))\\
&=& (x,y),
\end{eqnarray*}
}
and (3) is proved.  
\end{proof}
The advantage of the \wdl\ $L$ constructed in Lemma~\ref{L:combinehk} is that, in addition to extending the weakly and dual weakly complemented lattice structures induced by $h$  and $k$, it also keeps track of the closure and kernel systems.
\begin{definition}
Let $L$ be a bounded lattice and $x\in L$. The element $x^*\in L$ (resp. $x^+\in L$) is the \neu{pseudocomplement} (resp. \neu{dual pseudocomplement}) of $x$ if  
\[x\wedge y=0\iff y\leq x^*\quad\text{ (resp. $x\vee y=1\iff y\geq x^+$) for all $y\in L$}.\]
A \neu{double p-algebra} is a lattice in which every element has a pseudocomplement and a dual pseudocomplement. 
\end{definition}
\begin{example}
Boolean algebras are double p-algebras. Finite distributive lattices are double {p-algebras}. All distributive double p-algebras are weakly dicomplemented lattices.  $N_5$ is a double p-algebra that is not distributive. The double p-algebra operation $(^+,^*)$ on $N_5$ is however not a weak dicomplementation.   
\end{example}
The following result give a class of ``more concrete'' {\wdl}s, and can serve as prelude to the representation problem for {\wdl}s.
\begin{proposition}\label{L:JML}
Let $L$ be a finite lattice. Denote by $J(L)$ the set of join irreducible elements of $L$ and by $M(L)$ the set of meet irreducible elements of $L$ respectively. Define two unary operations $\tu$ and $\td$ on $L$ by 
$$x\tu:=\bigvee\{a\in J(L)\mid a\nleq x\}\quad \text{ and }\quad x\td:=\bigwedge\{m\in M(L)\mid m\ngeq x\}.$$
Then $(L,\wedge,\vee,\tu,\td,0,1)$ is a weakly dicomplemented lattice. In general, for $G\supseteq J(L)$ and $H\supseteq M(L)$, the operations $\sp{\tru\sb{G}}$ and $\sp{\trd\sb{H}}$ defined by
\[x\sp{\tru_G}:=\bigvee\{a\in G\mid a\nleq x\}\quad \text{ and }\quad x\sp{\trd_H}:=\bigwedge\{m\in H\mid m\ngeq x\}\]
turn $(L,\wedge,\vee,\sp{\tru_G},\sp{\trd_H},0,1)$ into a weakly dicomplemented lattice.
\end{proposition}
\begin{proof}
Let $G\supseteq J(L)$, $b\in G$ and $x\in L$. Then $b\nleq\bigvee\{a\in G\mid a\nleq x\}$ implies $b\leq x$; \newline i.e., $b\nleq x\sp{\tru_G}\implies b\leq x$. 
Thus $x\sp{\tru_G\tru_G}=\bigvee\{b\in G\mid b\nleq x\sp{\tru_G}\}\leq x$ and $(1)$ is proved. For $x\leq y$ we have $\{a\in G\mid a\nleq x\}\supseteq\{a\in G\mid a\nleq y\}$ implying $x\sp{\tru_G}\geq y\sp{\tru_G}$, and (2) is proved. For (3), it is enough to prove that for $a\in J(L)$, $a\leq x\iff a\leq (x\wedge y)\vee(x\wedge y\sp{\tru_G})$, since $J(L)$ is $\bigvee$-dense in $L$. Let $a\leq x$. 
We have $a\leq y$ or $a\leq y\sp{\tru_G}$. Then $a\leq x\wedge y$ or $a\leq x\wedge y\sp{\tru_G}$. Thus $a\leq (x\wedge y)\vee(x\wedge y\sp{\tru_G})$. The reverse inequality is obvious. $(1')-(3')$ are proved similarly.  
\end{proof}
Example~\ref{L:JML} above is a special case of concept algebras. Before we introduce concept algebras, let us recall some basic notions from FCA. The reader is referred to \cite{GW99}. As we mentioned before, FCA is based on the formalization of the notion of concept and concept hierarchy. Traditional philosophers considered a \neu{concept} to be determined by its extent and its intent. The extent consists of all objects belonging to the concept while the intent is the set of all attributes shared by all objects of the concept. In general, it may be difficult to list all objects or attributes of a concept. Therefore a specific {\em context} should be fixed to enable formalization. 
A \neu{formal context} is a triple $(G,M,I)$ of sets such that $I\subseteq G\times M$. The members of $G$ are called \neu{objects} and those of $M$ \neu{attributes}. If $(g,m)\in I$, then the object $g$ is said to have $m$ as an attribute. For subsets $A\subseteq G$ and $B\subseteq M$, $A'$ and $B'$ are  defined by
\[A':=\{m\in M\mid\ \forall g\in A\ g\I m \}\quad\text{and}\quad B':=\{g\in G\mid\ \forall m\in B\ g\I m \}.\]
A \neu{formal concept} of the formal context $(G,M,I)$ is a pair $(A,B)$ with $A\subseteq G$ and $B\subseteq M$ such that $A'=B$ and $B'=A$. The set $A$ is called the \neu{extent} and $B$ the \neu{intent} of the concept $(A,B)$. \neu{$\mathfrak{B}(G,M,I)$} denotes the set of all formal concepts of the formal context $(G,M,I)$. The concept hierarchy states that a concept is more general if it contains more objects. For capturing this notion a \neu{subconcept-superconcept relation} is defined: a concept $(A,B)$ is called a \neu{subconcept} of a concept $(C,D)$ provided that $A\subseteq C$ (which is equivalent to $D\subseteq B$); in this case, $(C,D)$ is a called \neu{superconcept} of $(A,B)$ and we write $(A,B)\leq(C,D)$. Obviously the subconcept-superconcept relation is an order relation on the set $\frak{B}(G,M,I)$ of all concepts of the formal context $(G,M,I)$.
The following result describing the concept hierarchy is considered as the basic theorem of FCA. 
\begin{theorem}[\cite{Wi82}]\label{T:BasicThmFCA} The poset $(\frak{B}(G,M,I),\leq)$ is a complete lattice in which infimum and supremum are given by:
{\small \[
\bigwedge_{t\in T}\left(A_t,B_t\right)=\left(\bigcap_{t\in
  T}A_t,\left(\bigcup_{t\in T}B_t\right)''\right) \text{ and }
\bigvee_{t\in T}\left(A_t,B_t\right)=\left(\left(\bigcup_{t\in
  T}A_t\right)'',\bigcap_{t\in T}B_t\right).
\] }
A complete lattice $L$ is isomorphic to $\frak{B}(G,M,I)$ iff
there are mappings $\tilde{\gamma}:G\to L$ and $\tilde{\mu}:M\to L$
such that $\tilde{\gamma}(G)$ is supremum-dense, $\tilde{\mu}(M)$ is
infimum-dense and $g\I m\iff \tilde{\gamma}g\leq\tilde{\mu}m$ for all
$(g,m)\in G\times M$. 
\end{theorem}
The poset $(\mathfrak{B}(G,M,I);\leq )$ is called the \neu{concept lattice} of the context $(G,M,I)$ and is denoted by $\underline{\frak{B}}(G,M,I)$. By Theorem~\ref{T:BasicThmFCA}, all complete lattices are (copies of) concept lattices. We adopt the notations below for $g\in G$ and $m\in M$:
\[g':=\{g\}',\quad m':=\{m\}',\quad {\gamma}g:=(g'',g')\quad\mbox{ and }\quad {\mu}m:=(m',m'').\]
The concept $\gamma g$ is called \neu{object concept} and $\mu m$ \neu{attribute concept}. They form the building blocks of the concept lattice. The sets $\gamma G$ is supremum-dense and $\mu M$ is infimum-dense in $\frak{B}(G,M,I)$. We usually assume our context clarified, meaning that
 $$x'=y'\implies x=y \quad\text{ for all } x,y \text{ in } G\cup M.$$
If $\gamma g$ is supremum-irreducible we say that the \neu{object} $g$ is \neu{irreducible}. An \neu{attribute} $m$ is called \neu{irreducible} if the attribute concept $\mu m$ is infimum-irreducible. A \neu{formal context} is called \neu{reduced} if all its objects and attributes are irreducible. For every finite nonempty lattice $L$ there is, up to isomorphism, a unique reduced context ${\KK(L):=(J(L),M(L),\leq)}$ such that $L\cong\underline{\frak{B}}(\KK(L))$. We call it \neu{standard context} of $L$. The meet and join operations in the concept lattice can be used to formalize respectively the conjunction and disjunction on concepts~\cite{GW99a}. To formalize the negation, the main problem is that the complement of an extent is probably not and extent and the complement of an intent might not be an intent. Therefore two operations are introduced as follows:
\begin{definition}\label{D:ca}
Let $\KK:=(G,M,I)$ be a formal context. We define for each formal concept $(A,B)$ 
   \[\mbox{ its \neu{weak negation} by}\qquad (A,B)\tu :=\left(\left(G\setminus A\right)'',\left(G\setminus A\right)'\right)\]
   \[\mbox{and its \neu{weak opposition} by}\quad (A,B)\td :=\left(\left(M\setminus B\right)',\left(M\setminus B\right)''\right).\]
   $\frak{A}(\KK):=\left(\frak{B}(\KK);\wedge,\vee,\tu,\td,0,1\right)$ is called the \neu{concept algebra} of the formal context $\KK$, where $\wedge$ and $\vee$ denote the meet and the join operations of the concept lattice.
   \end{definition}
\noindent These operations satisfy the equations in Definition~\ref{D:wdl} (cf.~\cite{Wi00}). In fact, concept algebras are typical examples of weakly dicomplemented lattices. One of the important and still unsolved problems in this topic is to find out the equational theory of concept algebras; that is the set of all equations valid in all concept algebras. Is it finitely generated? I.e. is there a finite set $\mathcal{E}$ of equations valid in all concept algebras such that each equation valid in all concept algebras follows from $\mathcal{E}$? We start with the set of equations defining a \wdl\ and have to check whether they are enough to represent the equational theory of concept algebras. 
 This problem, known as ``representation problem'' (\cite{Kw04}), can be split in three sub-problems:
\begin{itemize}
\item[SRP] \neu{Strong representation problem}: describe weakly dicomplemented lattices that are isomorphic to concept algebras.
\item[EAP] \neu{Equational axiomatization problem}: find a set of equations that generate the equational theory of concept algebras.
\item[CEP] \neu{Concrete embedding problem}: given a weakly dicomplemented lattice $L$, is there a context $\KK\sb{\trd}\sp{\tru}(L)$ such that $L$ embeds into the concept algebra of $\KK\sb{\trd}\sp{\tru}(L)$?
\end{itemize}
We proved (see \cite{Kw04} or \cite{GK07}) that finite distributive weakly dicomplemented lattices are isomorphic to concept algebras. However we cannot expect all weakly dicomplemented lattices to be isomorphic to concept algebras, since concept algebras are necessary complete lattices. In Section~\ref{S:main} we will show that even being complete is not enough for weakly dicomplemented lattices to be isomorphic to concept algebras. Before that we show in Section~\ref{S:wdl+neg} that {\wdl}s generalize Boolean algebras. 
\section{Weakly Dicomplemented Lattices with Negation}\label{S:wdl+neg}
Example~\ref{E:ba} states that duplicating the
complementation of a Boolean algebra leads to a \wdl. Does the
converse hold? I.e., is a \wdl\ in which the weak complementation and the dual weak complementation are duplicate a Boolean algebra? The finite case is easily
obtained~[Co\-rollary~\ref{C:fwdl+neg}]. Major parts of this section are taken from \cite{Kw04}. We will also describe {\wdl}s whose Boolean part is the intersection of their skeletons (definitions below). 
\begin{definition}\label{D:wdl+neg}
   A weakly dicomplemented lattice is said to be \neu{with
   negation}\index{lattice!weakly dicomplemented!with negation}
   if the unary operations coincide, i.e., if $x\td = x\tu$ for all
   $x$. In this case we set $x\tu=:\bar{x}:=x\td$.
   \end{definition}
   \begin{lemma}\label{P:wdl+neg-uc}
   A \wdl\ with negation is uniquely complemented.
   \end{lemma}
   \begin{proof}
$x^{\tru\tru}\leq x\leq x^{\trd\trd}$ implies that $x=\bar{\bar{x}}$. Moreover,
   $x\wedge \bar{x} =0$ and $\bar{x}$ is a complement of $x$. If $y$ is another complement of $x$ then
   \[ x=(x\wedge y)\vee(x\wedge\bar{y})= x\wedge\bar{y} \implies x\leq \bar{y}\]
   \[ x=(x\vee y)\wedge(x\vee \bar{y})= x\vee \bar{y} \implies x\geq \bar{y}\]
   Then $\bar{y} = x$ and $\bar{x} = y$. $L$ is therefore a uniquely complemented lattice.  
   \end{proof}
It can be easily seen that each uniquely complemented atomic lattice
is a copy of the power set of the set of its atoms, and therefore
distributive. Thus 
\begin{corollary}\label{C:fwdl+neg}
 The finite weakly dicomplemented lattices with negation are exactly
 the finite Boolean algebras. 
\end{corollary}
Of course, the natural question will be if the converse of
Lemma~\ref{P:wdl+neg-uc} holds. That is, can any uniquely
complemented lattice be endowed with a structure of a weakly
dicomplemented lattice with negation? The answer is yes for
distributive lattices. If the assertion of Corollary~\ref{C:fwdl+neg}
can be extended to lattices in general, the answer will
unfortunately be no. In fact R.~P.~Dilworth proved that each lattice can be embedded into a uniquely complemented lattice \cite{Di45}. The immediate consequence is the existence of non-distributive uniquely complemented lattices. They are however infinite. 
If a uniquely complemented lattice could be endowed with a structure of weakly
dicomplemented lattice, it would be distributive. This cannot be true
for non distributive uniquely complemented lattices. 

\begin{lemma}\label{L:wdl+neg-dml}
Each weakly dicomplemented lattice with negation $L$
satisfies the  de Morgan laws. 
\end{lemma}
\begin{proof}
We want to prove that $\overline{x\wedge y}=\bar{x}\vee\bar{y}$.
\[
(\bar{x}\vee\bar{y})\vee(x\wedge y)\geq \bar{x}\vee(x\wedge \bar{y})\vee(x\wedge
y)=\bar{x}\vee x=1 
\]
 and
\[
(\bar{x}\vee\bar{y})\wedge(x\wedge y)\leq (\bar{x}\vee \bar{y})\wedge x\wedge(\bar{x}\vee y)=\bar{x}\wedge x=0.
\]
So $\bar{x}\vee \bar{y}$ is a complement of $x\wedge y$, hence by uniqueness
it is equal to $\overline{x\wedge y}$. Dually we have $\overline{x\vee y}=\bar{x}\wedge \bar{y}$.  
\end{proof}
Now for the distributivity we can show that
\begin{lemma}\label{P:wdl+neg-distr}
$\overline{x\wedge(y\vee z)}$ is a complement of $(x\wedge y)\vee(x\wedge z)$.
\end{lemma}
\begin{proof}
Since in every lattice the equation $x\wedge(y\vee z)\geq(x\wedge y)\vee(x\wedge z)$ holds, we have that $\overline{x\wedge(y\vee z)}\leq\overline{(x\wedge y)\vee(x\wedge z)}$; so we have only to show that 
$$\overline{x\wedge(y\vee z)}\vee(x\wedge y)\vee(x\wedge z)=1.$$
Using the de Morgan laws and axiom (3) several times we obtain:
\begin{eqnarray*}
\overline{x\wedge(y\vee z)}\vee(x\wedge y)\vee(x\wedge z)
&=&\bar{x}\vee(\bar{y}\wedge \bar{z})\vee(x\wedge y)\vee(x\wedge z)\\
&=& \bar{x}\vee(\bar{y}\wedge \bar{z}\wedge x)\vee(\bar{y}\wedge \bar{z}\wedge
\bar{x})\vee(x\wedge y\wedge z)\\
& &\vee(x\wedge y\wedge \bar{z}) \vee(x\wedge z\wedge \bar{y})\\
&=& \bar{x}\vee(\bar{y}\wedge \bar{z}\wedge \bar{x})\vee(x\wedge y\wedge z)\vee(x\wedge y\wedge \bar{z})\\
& &\vee(x\wedge \bar{y}\wedge z)\vee(x\wedge \bar{y}\wedge \bar{z})\\
&=& \bar{x}\vee(\bar{y}\wedge \bar{z}\wedge \bar{x})\vee(x\wedge y)\vee(x\wedge \bar{y})\\
&=& \bar{x}\vee(\bar{y}\wedge \bar{z}\wedge \bar{x})\vee x\\
&=& 1.
\end{eqnarray*}
Thus $\overline{x\wedge(y\vee z)}$ is a complement of $(x\wedge y)\vee(x\wedge z)$.  
\end{proof}
Since the complement is unique we get the
equality 
\[
x\wedge(y\vee z)=\overline{\overline{x\wedge(y\vee z)}}=(x\wedge y)\vee(x\wedge z).
\]
Thus {\wdl}s generalize Boolean algebras in the following sense
\begin{theorem}\label{T:wdl+neg}
Boolean algebras with duplicated complementation\footnote{see Example~\ref{E:ba}} are {\wdl}s. If $\tu=\td$ in a \wdl\ $L$ then $(L,\wedge,\vee,\tu,0,1)$ is a Boolean algebra. 
\end{theorem}
As the equality $x\tu=x\td$ not always holds, we can look for maximal substructures with this property.
\begin{definition} For any weakly dicomplemented lattice $L$, we will call $$B(L):=\{x\in L\mid x\tu=x\td\}$$ the \neu{subset of elements with negation}.
\end{definition}
As in Definition~\ref{D:wdl+neg} we denote by $\bar{x}$ the common value of
$x\tu$ and $x\td$, for any $x\in B(L)$. We set 
$$L\tu:=\{a\tu\mid a\in L\}=\{a\in L\mid a{\tu}{\tu}=a\}$$
 and call it the \neu{skeleton} of $L$, as well as 
$$L\td:=\{a\td\mid a\in L\}=\{a\in L\mid a{\td}{\td}=a\}$$
 and call it the \neu{dual skeleton} of $L$.
\begin{corollary}\label{C:elts+neg} $(B(L),\wedge,\vee,\bar{\phantom{a}},0,1)$ is a Boolean algebra that is a subalgebra of the skeleton and the dual skeleton.
\end{corollary}
\begin{proof}
From $x\tu=x\td$ we get $x\sp{\tru\tru}=x\sp{\trd\tru}$ and
  $x\sp{\tru\trd}=x\sp{\trd\trd}$. Thus
  \[
x\sp{\tru\trd}=x\sp{\tru\tru}=x=x\sp{\trd\trd}=x\sp{\trd\tru}
\]
 and
  $B(L)$ is closed under the operations $\tu$ and $\td$. We will prove
  that $B(L)$ is a subalgebra of $L$. We consider $x$ and $y$ in
  $B(L)$.  
We have 
\[
(x\wedge y)\tu=x\tu\vee y\tu=x\td\vee y\td \leq (x\wedge y)\td
\leq(x\wedge y)\tu \mbox{ and }
\]
\[
(x\vee y)\td=x\td\wedge y\td=x\tu\wedge y\tu \geq (x\vee y)\tu \geq(x\vee y)\td.\]
Thus $x\wedge y$ and $x\vee y$ belong to $B(L)$. $B(L)$ is a weakly
dicomplemented lattice with negation, and is by
Theorem~\ref{T:wdl+neg}, a Boolean algebra.  
\end{proof} 
While proving Corollary~\ref{C:elts+neg} we have also shown that $B(L)$ is a
subalgebra of $L$. It is, in fact, the largest Boolean algebra that is a
subalgebra of the skeletons and of $L$. We call it the
\neu{Boolean part} of $L$. The inclusion $B(L)\subseteq L\tu\cap L\td$ can be strict. For the \wdl\ $L_1$ in Fig.~\ref{Fi:liebsp}, we have 
\[B(L_1)=\{0,1\},\ L_1\tu=\{0,1,c,d,e,c\tu,d\tu,e\tu\} \text{ and } L_1\td=\{0,1,c,a,b,c\td,a\td,b\td\}.\] 
 Thus $B(L_1)\subsetneq L_1\tu\cap L_1\td$. It would be nice to find under which conditions the Boolean part is the intersection of the skeleton and dual skeleton. 
\begin{lemma}
If $L$ is a finite distributive lattice with $\td=\ast$ (pseudocomplementation) and $\tu=+$ (dual pseudocomplementation), then $B(L)$ is the set of complemented elements of $L$.
\end{lemma}
\begin{proof}
 Let $L$ be a finite distributive lattice with $\td=\ast$ and $\tu=+$. We denote by $C(L)$ the set of complemented elements of $L$. Of course $B(L)\subseteq C(L)$. Let $x\in C(L)$. From the distributi\-vi\-ty there is a unique elements $z\in L$ such that $x\vee z=1$ and $x\wedge z=0$. Then $z\leq x\td\leq x\tu\leq z$, and $x\in B(L)$.  
\end{proof}
Even in this case, the Boolean part can still be strictly smaller than the intersection of the skeletons. For $L_1$ in Fig.~\ref{Fi:liebsp} we have
$$B(L_1)\subsetneq L_1\tu\cap L_1\td=\{0,1,c,a\td\}= C(L_1).$$ 
For $L_2$ in Fig~\ref{Fi:liebsp}, we have $\tu=\sp{+}$ and $\td=\sp{\ast}$; but \[L_2\tu=\{0,1,c,c\tu\},\ L_2\td=\{0,1,c,c\td\},\ B(L_2)=\{0,1\}=C(L_2)\subsetneq\{0,1,c\}=L_2\tu\cap L_2\td.\] 

\begin{figure}
\begin{minipage}[t]{6cm}
\begin{center}
\begin{texdraw}
\drawdim{cm} \setunitscale .6 
\move(4 2) \lvec(1 5)\lvec(4 8) \lvec(7 5)\lvec(4 2)
\move(5 3) \lvec(2 6) 
\move(6 4) \lvec(3 7) 
\move(3 3) \lvec(6 6)
\move(2 4) \lvec(5 7) 
\move(4 2) \fcir f:0 r:0.2 \lcir r:0.2
\move(3 3) \fcir f:0.7 r:0.2 \lcir r:0.2 
\move(2 4) \fcir f:0.7 r:0.2 \lcir r:0.2 
\move(1 5) \fcir f:0.7 r:0.2 \lcir r:0.2 
\move(5 3) \fcir f:0.7 r:0.2 \lcir r:0.2 
\move(4 4) \fcir f:1 r:0.2 \lcir r:0.2 
\move(3 5) \fcir f:1 r:0.2 \lcir r:0.2 
\move(2 6) \fcir f:0.7 r:0.2 \lcir r:0.2 
\move(6 4) \fcir f:0.7 r:0.2 \lcir r:0.2
\move(5 5) \fcir f:1 r:0.2 \lcir r:0.2 
\move(4 6) \fcir f:1 r:0.2 \lcir r:0.2 
\move(3 7) \fcir f:0.7 r:0.2 \lcir r:0.2 
\move(7 5) \fcir f:0.7 r:0.2 \lcir r:0.2 
\move(6 6) \fcir f:0.7 r:0.2 \lcir r:0.2 
\move(5 7) \fcir f:0.7 r:0.2 \lcir r:0.2 
\move(4 8) \fcir f:0 r:0.2 \lcir r:0.2

\htext(4.3 2){$0$} \htext(4.3 8){$1$} \htext(2.5 3){$a$}
\htext(1.5 4){$b$} \htext(0.5 5){$c$} \htext(1.5 6){$d$}
\htext(2.5 7){$e$} \htext(5.3 3){$c\td$} \htext(6.3 4){$b\td$}
\htext(7.3 5){$a\td=e\tu$} \htext(6.3 6){$d\tu$} \htext(5.3
7){$c\tu$} \htext(4.3 4){$u$} \htext(5.3 5){$v$} \htext(3.3 5){$w$}
\htext(4.3 6){$z$}
\htext(2.3 2){$L_1$}
\end{texdraw}
\end{center}
\end{minipage}
\begin{minipage}[t]{6cm}
\begin{center}
\begin{texdraw}
\drawdim{cm} \setunitscale .6
\move(4 2) \lvec(1 5)\lvec(4 8) \lvec(6 6)
\move(6 4) \lvec(4 2)
\move(5 3) \lvec(2 6) 
\move(6 4) \lvec(3 7) 
\move(3 3) \lvec(6 6)
\move(2 4) \lvec(5 7) 
\move(4 2) \fcir f:0 r:0.2 \lcir r:0.2
\move(3 3) \fcir f:0.7 r:0.2 \lcir r:0.2 
\move(2 4) \fcir f:0.7 r:0.2 \lcir r:0.2 
\move(1 5) \fcir f:0.7 r:0.2 \lcir r:0.2 
\move(5 3) \fcir f:0.7 r:0.2 \lcir r:0.2 
\move(4 4) \fcir f:1 r:0.2 \lcir r:0.2 
\move(3 5) \fcir f:1 r:0.2 \lcir r:0.2 
\move(2 6) \fcir f:0.7 r:0.2 \lcir r:0.2 
\move(6 4) \fcir f:0.7 r:0.2 \lcir r:0.2
\move(5 5) \fcir f:1 r:0.2 \lcir r:0.2 
\move(4 6) \fcir f:1 r:0.2 \lcir r:0.2 
\move(3 7) \fcir f:0.7 r:0.2 \lcir r:0.2 
\move(6 6) \fcir f:0.7 r:0.2 \lcir r:0.2 
\move(5 7) \fcir f:0.7 r:0.2 \lcir r:0.2 
\move(4 8) \fcir f:0 r:0.2 \lcir r:0.2

\htext(4.3 2){$0$} \htext(4.3 8){$1$} \htext(2.5 3){$a$}
\htext(1.5 4){$b$} \htext(0.5 5){$c$} \htext(1.5 6){$d$}
\htext(2.5 7){$e$} \htext(5.3 2.8){$f$} \htext(6.3 4){$c\td$}
\htext(6.3 6){$c\tu$} \htext(5.3 7){$g$} \htext(4.3 4){$u$} \htext(5.3 5){$v$} \htext(3.3 5){$w$} \htext(4.3 6){$z$}
\htext(2.3 2){$L_2$}
\end{texdraw}
\end{center}
\end{minipage}
\caption{Examples of dicomplementations. For $L_1$, the elements $c$, $b$ and $a$ are each image (of their image). The operation $\tu$ is the dual of $\td$. 
For $L_2$, $\tu=\sp{+}$ and $\td=\sp{\ast}$.}\label{Fi:liebsp} 
\end{figure}
\begin{lemma}
$B(L)=L\tu\cap L\td$ iff $x{\tu}{\tu}=x{\td}{\td}\implies x{\tu}{\td}=x{\td}{\tu}$. 
\end{lemma}
\begin{proof}
($\Rightarrow$). Let $x\in L$ such that $x{\tu}{\tu}=x{\td}{\td}$. Then $x\in L\tu\cap L\td=B(L)$ and implies $x\tu=x\td$. Therefore 
$x{\tu}{\td}=x{\td}{\td}=x=x{\tu}{\tu}=x{\td}{\tu}$.

($\Leftarrow$). Let $x\in L\tu\cap L\td$. Then $x{\tu}{\tu}=x=x{\td}{\td}$ and implies $x{\tu}=x{\td}{\tu}{\tu}\leq x\td$. Thus $x\tu=x\td$, and $x\in B(L)$.  
\end{proof}
\begin{lemma}
If $L\tu$ and $L\td$ are subalgebras of $L$, then there are complemented.  
\end{lemma}
\begin{proof}
We assume that $L\td$ is a subalgebra of $L$. Let $x\in L\td$. Then $x\wedge x\td=0$ and $x\vee x\td=t\td$ for some $t\in L$. Therefore 
\[0=(x\vee x\td)\td=t^{\trd\trd} \implies 1=0\td=t^{\trd\trd\trd}=t\td=x\vee x\td.\]
Thus $L\td$ is complemented. The proof for $L\tu$ is obtained analogously.  
\end{proof}
In general, $L\tu$ and $L\td$ are orthocomplemented lattice, when considered as lattice on their own~\cite{Kw04}.  
\section{Strong representation problem}\label{S:main}
We start this section by a negative result, namely by showing that {\em completeness is not enough for weakly dicomplemented lattices to be (copies of) concept algebras}. 
\begin{theorem}\label{T:main}
There is no formal context whose concept algebra is isomorphic to a complete atomfree Boolean algebra.
\end{theorem}
\begin{proof}
 Let $B$ be a complete and atomfree Boolean algebra. By Theorem~\ref{T:BasicThmFCA}, there is a context $(G,M,I)$ such that $\underline{\frak{B}}(G,M,I)\cong B$ (lattice isomorphism). Without loss of generality, we can assume that $(G,M,I)$ is a subcontext of $(B,B,\leq)$. We claim that there are $g,h\in G$ with $0<h<g<1$. In fact, for an element $g\in G\subseteq B$ with $0\neq g$ there is $a\in B$ such that $0<a<g$, since $B$ is atomfree. Moreover $G$ is $\bigvee$-dense in $B$ and then $0\neq a=\bigvee\{x\in G\mid x\leq a\}$, implying that $\{x\in G\mid 0<x\leq a\}\neq\emptyset$. Thus we can choose $h\in G$ with $0<h\leq a<g$. In the concept algebra of $(G,M,\leq)$ we have $h\tu=\bigvee\{x\in G\mid x\nleq h\}\geq g>h$. From $h\vee h\tu=1$ we get $h\tu=1\neq \bar{h}$ (the complement of $h$ in $B$).  
\end{proof}
Theorem~\ref{T:main} says that an atomfree Boolean algebra is not isomorphic to a concept algebra. However it can be embedded into a concept algebra. The corresponding context is constructed via ultrafilters. A general construction was presented in \cite{Kw04}.

\def\Fpr{\frak{F}\sb{\rm pr}}
\def\Ipr{\frak{I}\sb{\rm pr}}
\def\Fx#1{\cal{F}\sb{\rm #1}}
\def\Ix#1{\cal{I}\sb{\rm #1}}

\begin{definition}\label{D:primaryfi}
A \neu{primary filter} is a (lattice) filter that contains $w$ or $w\tu$ for all $w\in L$. Dually, a \neu{primary ideal} is an ideal that contains $w$ or $w\td$ for all $w\in L$. $\Fpr(L)$ denotes the set of all primary filters and $\Ipr(L)$ the set of primary ideals  of $L$.
\end{definition}
For Boolean algebras, a proper filter $F$ is primary iff it is an ultrafilter, iff it is a prime filter ($x\vee y\in F \implies x\in F$ or $y\in F$). The following result based on Zorn's lemma provides the sets of a context $\KK\sb{\trd}\sp{\tru}(L)$ which is the best candidate for representing a \wdl\ $L$.
\begin{theorem}[``Prime ideal theorem''~\cite{Kw04}]\label{T:pit}
For every filter $F$ and every ideal $I$ such that $F\cap I=\emptyset$
there is a primary filter $G$ containing $F$ and disjoint from $I$. Dually, for every ideal $I$ and every filter $F$ such that $I\cap F=\emptyset$
there is a primary ideal $J$ containing $I$ and disjoint from $F$.
\end{theorem}
\begin{corollary}\label{C:pit}
If $x\not\le y$ in $L$, then there exists a primary filter $F$
containing $x$ and not $y$. 
\end{corollary}
For $x\in L$, we set
\[\Fx{x}:=\{F\in\Fpr(L)\mid x\in F\}\quad\text{ and }\quad
\Ix{x}:=\{I\in\Ipr(L)\mid x\in I\}.\]
The \neu{canonical context} of a weakly dicomplemented lattice $L$ is the formal context  
\[\bb{K}\sp{\tu}\sb{\td}(L):=(\Fpr(L),\Ipr(L),\mathrel{\Box})\quad\text{ with }F\mathrel{\Box}I:\iff F\cap I\ne\emptyset.\]
The derivation in $\bb{K}\sp{\tu}\sb{\td}(L)$ yields, $\Fx{x}'=\Ix{x}$ and $\Ix{x}'=\Fx{x}$ for every $x\in L$. Moreover, the map 
\begin{eqnarray*}
 i \colon L & \to & \frak{B}\left(\bb{K}\sp{\tu}\sb{\td}(L)\right) \\ 
    x & \mapsto  & (\Fx{x},\Ix{x})
\end{eqnarray*}
is a bounded lattice embedding with 
\[
i(x\td)\leq i(x)\td\leq i(x)\tu \leq i(x\tu).
\]
 If the first and last inequalities above were equalities, we would get a weakly dicomplemented lattice embedding into the concept algebra of $\KK\sp{\tru}\sb{\trd}(L)$. This would be a solution to the representation problem of weakly dicomplemented lattices.
\begin{theorem}\label{T:atcompletion}
If $L$ is a Boolean algebra, then the concept algebra of $\bb{K}_\trd^\tru(L)$ is a complete and atomic Boolean algebra into which $L$ embeds.
\end{theorem}
\begin{proof}
If $B$ is a Boolean algebra, then a proper filter $F$ of $L$ is primary iff it is an ultrafilter, and a proper ideal $J$ is primary iff it is maximal. Thus $\Fpr(L)$ is the set of ultrafilters of $L$ and $\Ipr(L)$ the set of its maximal ideals. In addition, the complement of an ultrafilter is a maximal ideal and vice-versa. For $F\in\Fpr(L)$, $L\setminus F$ is the only primary ideal that does not intersect $F$, and for any $J\in\Ipr(L)$, $L\setminus J$ is the only primary filter that does not intersect $J$. Thus the context  $\bb{K}_\trd^\tru(L)$ is a copy of $(\Fpr(L),\Fpr(L),\neq)$. The concepts of this context are exactly pairs $(A,B)$ such that $A\cup B=\Fpr(L)$ and $A\cap B=\emptyset$. Thus $\frak{B}(\bb{K}_\trd^\tru(L))\cong\cal{P}(\Fpr(L))$ and each subset $A$ of $\Fpr(L)$ is an extent of $\bb{K}_\trd^\tru(L)$. It remains to prove that the lattice embedding 
\begin{eqnarray*}
 i \colon L & \to & \frak{B}\left(\bb{K}\sp{\tu}\sb{\td}(L)\right) \\ 
    x & \mapsto  & (\Fx{x},\Ix{x})
\end{eqnarray*}
is also a Boolean algebra embedding. If $i(x\tu)\neq i(x)\tu$ then there is 
\[
F\in\Fx{x\tu}\setminus\left(\Fpr(L)\setminus\Fx{x}\right)''=\Fx{x\tu}\setminus\left(\Fpr(L)\setminus\Fx{x}\right)=\emptyset,
\]
 which is a contradiction. Similarly $i(x\td)=i(x)\td$. Therefore $B$ embeds into the complete and atomic Boolean algebra $\frak{A}\left(\bb{K}_\trd^\tru(L)\right)$ which is a copy of $\mathcal{P}\left(\Fpr(L)\right)$.  
\end{proof}
The above result is a new proof to a well-known result (Corollary~\ref{C:Stone}) due to Marshall Stone~\cite{St36}. The advantage here is that the proof is simple and does not require any know\-ledge from topology. Recall that a field of subsets of a set $X$ is a subalgebra of $\cal{P}(X)$, i.e. a family of subsets of $X$ that contains $\emptyset$ and $X$, and that is closed under union, intersection, and complementation.
\begin{corollary}[\cite{St36}]\label{C:Stone}
Each Boolean algebra embeds into a field of sets. 
\end{corollary}
We conclude this section by an example.
Consider the Boolean algebra ${F}\NN$ of finite and cofinite subsets of $\NN$. It is not complete. But $\cal{P}(\NN)$ is a complete and atomic Boolean algebra containing ${F}\NN$. 
 By Theorem~\ref{T:atcompletion} $\frak{A}(\KK_\trd^\tru({F}\NN))$ is also a complete and atomic Boolean algebra into which ${F}\NN$ embeds. The atoms of ${F}\NN$ are $\{n\}, n\in\NN$. These generate its principal ultrafilters. ${F}\NN$ has exactly one non-principal ultrafilter $U$ (the cofinite subsets). Thus $|{F}\NN|=|\NN|+1=|\NN|$. We can find a bijection let say $f$ between the atoms of $\cal{P}(\NN)$ and the atoms of $\frak{A}(\KK_\trd^\tru({F}\NN))$. $f$ induces an isomorphism $\hat{f}:\cal{P}(\NN)\to\frak{A}(\KK_\trd^\tru({F}\NN))$. Henceforth, it is natural to look for a universal property to characterize $\frak{A}(\KK_\trd^\tru(B))$ for any Boolean algebra $B$. For example {\em is $\frak{A}(\KK_\trd^\tru(B))$ the smallest complete and atomic Boolean algebra into which $B$ embeds}?
\section{Conclusion}
Weakly dicomplemented lattices with negation are exactly Boolean algebras (Theorem~\ref{T:wdl+neg}). Even if they are not always isomorphic to concept algebras (Theorem~\ref{T:main}), they embed into concept algebras (Theorem~\ref{T:atcompletion}). Finite distributive {\wdl}s are isomorphic to concept algebras~\cite{GK07}. Extending these results to finite {\wdl}s in one sense and to distributive {\wdl}s in the other are the next steps towards the representation of weakly dicomplemented lattices. Finding a kind of universal property to characterize the construction in Theorem~\ref{T:atcompletion} is a natural question to be addressed. 
\section*{Acknowledgments}
This work was initiated as the first author was at TU Dresden, and was carried out while the first author was visiting Hitotsubashi University, Kunitachi, Tokyo, Japan, supported by the grant of the second author. We would like to thank both  universities for providing working facilities.
The authors would like to thank P\'eter P\'al P\'alfy for helpful discussions on the topic.


\begin{thebibliography}{00}
\bibitem[Bo54]{Bo54} G.\,Boole. \newblock \textsl{An investigation of the laws of thought on which are founded the ma\-the\-matical theories of logic and probabilities}. \newblock Macmillan 1854. Reprinted by Dover Publications, New York (1958).
\bibitem[BW06]{BW06} C.\,Burgmann \& R.\,Wille. \newblock \textsl{The Basic Theorem on Preconcept Lattices}. In R.\,Missaoui and J.\,Schmid (Eds.), \textsl{Formal Concept Analysis, 4th International Conference, ICFCA 2006, Dresden, Germany, February 13-17, 2006, Proceedings}. LNAI \textbf{3874}. \newblock Springer Verlag Berlin Heidelberg (2006) 80-88.
\bibitem[BH05]{BH05} P.\,Burmeister and R.\,Holzer. \newblock \textsl{Treating Incomplete Knowledge in Formal Concept Analysis}. In B.\,Ganter, G.\,Stumme and R.\,Wille (Eds.) \textsl{Formal Concept Analysis, Foundations and Applications}. LNAI~\textbf{3626} Springer (2005) 114-126.
\bibitem[BS81]{BS81} S.\,Burris \& H.\,P.\,Sankappanavar. \newblock \textsl{A course in universal algebra}. \newblock Springer Verlag (1981).
\bibitem[DP02]{DP02} B.\,A.\,Davey \& H.\,A.\,Priestley. \newblock \textsl{Introduction to lattices and order}. \newblock Second edition \newblock  Cambridge (2002).
\bibitem[DE98]{DE98} K.\,Deiters \& M.\,Ern\'e. \newblock \textsl{Negations and contrapositions of complete lattices}. \newblock  Discrete Mathematics~\textbf{181} No.~1-3 (1998) 91-111.
\bibitem[Di45]{Di45} R.\,P.\,Dilworth. \newblock \textsl{Lattices with unique complements}. \newblock Transactions of the American Mathematical Society~\textbf{57} (1945) 123-154.
\bibitem[Fe06]{Fe06} S.\,Ferr\'e. \newblock \textsl{Negation, Opposition, and Possibility in Logical Concept Analysis}. In R.\,Missaoui and J.\,Schmid (Eds.) \textsl{Formal Concept Analysis, 4th International Conference, ICFCA 2006, Dresden, Germany, February 13-17, 2006, Proceedings}. LNAI~\textbf{3874} Springer (2006) 130-145.
\bibitem[GD86]{GD86} J.\,L.\,Guigues and V.\,Duquenne. \newblock\textsl{Familles minimales d'implications informatives r\'esultant d'un tableau de donn\'ees binaires}. Math\'ematiques et Sciences Humaines~\textbf{95} (1986) 5-18.
\bibitem[GK07]{GK07} B.\,Ganter \& L.\,Kwuida. \newblock \textsl{Finite distributive concept algebras}. \newblock Order. (2007).
\bibitem[GW99]{GW99} B.\,Ganter \& R.\,Wille. \newblock \textsl{Formal Concept Analysis. Mathematical Foundations}. \newblock  Springer (1999).
\bibitem[GW99a]{GW99a} B.\,Ganter \& R.\,Wille. \newblock \textsl{Contextual attribute logic}. In W.\,Tepfenhart and W.\,Cyre (Eds.),\newblock \textsl{Conceptual Structures: Standards and Practices, 7th International Conference on Conceptual Structures, ICCS'99, Blacksburg, Virginia, USA, July  12-15, 1999, Proceedings}. LNAI~\textbf{1640}. Springer (1999) 337-338.
\bibitem[HLSW01]{HLSW} C.\,Herrmann, P.\,Luksch, M.\,Skorsky \& R.\,Wille. \newblock \textsl{Algebras of semiconcepts and double Boolean algebras}. \newblock Contributions to General Algebra~\textbf{13}. \textsl{Proceedings of the Dresden Conference 2000 and the Summer School 1999}. Verlag Johannes Heyn, Klagenfurt (2001).
\bibitem[Hl89]{Hl89} L.\,R.\,Horn. \newblock \textsl{A natural history of negation}. The University of Chicago Press, Chicago and London (1989).
\bibitem[Kw04]{Kw04} L.\,Kwuida. \newblock \textsl{Dicomplemented lattices. A contextual generalization of Boolean algebras}. \newblock Dissertation TU Dresden. Shaker Verlag. (2004).
\bibitem[Kw09]{Kw09} L.\, Kwuida. \newblock \textsl{Axiomatization of Boolean algebras via weak dicomplementations}. ArXiv:\textbf{0911.0200} (2009).  
\bibitem[KM08]{KM08} L.\,Kwuida \& H.\,Machida. \newblock \textsl{On the isomorphism problem of weakly dicomplemented lattices}. Proceedings of the sixth international conference on concept lattices and their applications CLA2008. Radim Belohlavek and Sergei O.\,Kuznetsov (Eds.) Palack\'y University, Olomouc, Czech Republic (2008) 217-228. \quad(ISBN: 978-80-244-2111-7)
\bibitem[MNR08]{MNR08} R.\,Missaoui, L.\,Nourine and Y.\,Renaud. \newblock \textsl{Generating Positive and Negative Exact Rules Using Formal Concept Analysis: Problems and Solutions}. In R.\,Medina and S.A.\,Obiedkov (Eds.) \textsl{Formal Concept Analysis, 6th International Conference, ICFCA 2008, Montreal, Canada, February 25-28, 2008, Proceedings}. LNAI~\textbf{4933} Springer (2008) 169-181.
\bibitem[Pr70]{Pr70} H.\,A.\,Priestley. \newblock \textsl{Representation of distributive lattices by means of ordered Stone spaces}. \newblock Bulletin of the London Mathematical Society~\textbf{2} (1970) 186-190. 
\bibitem[Pr06]{Pr06} U.\,Priss. \newblock \textsl{Formal Concept Analysis in Information Science}. In B.\,Cronin (Ed.) Annual Review of Information Science and Technology~\textbf{40} (2006) 521-543.
\bibitem[St36]{St36} M.\,H.\,Stone. \newblock \textsl{The theory of representations for Boolean algebras}. \newblock Transactions of the American Mathematical Society~\textbf{40} (1936) 37-111.
\bibitem[VW05]{VW05} B.\,Vormbrock \& R.\,Wille. \newblock \textsl{Semiconcept and Protoconcept Algebras: The Basic Theorems}. In B.\,Ganter et al. (Eds.), \textsl{Formal Concept Analysis, Foundations and Applications}, LNAI~\textbf{3626} \newblock Springer (2005) 34--48.
\bibitem[Wa96]{Wa96} H.\,Wansing (Ed.). \newblock \textsl{Negation: a notion in focus}. \newblock \textsl{Perspectives in analytical philosophy}. \textbf{7} de Gruyter (1996).
\bibitem[Wi82]{Wi82} R.\,Wille. \newblock \textsl{Restructuring lattice theory:  an approach based on hierarchies of concepts}. \newblock In I. Rival (Ed.) \newblock \textsl{Ordered Sets}. Reidel (1982) 445-470.
\bibitem[Wi96]{Wi96} R.\,Wille. \newblock \textsl{Restructuring mathematical logic: an approach based on Peirce's pragmatism}. In A.\, Ursini et al. (Eds.) \newblock \textsl{Logic and algebra}. Marcel Dekker. Lecture Notes in Pure and Applied Mathematics~\textbf{180} (1996) 267-281.
\bibitem[Wi00]{Wi00} R.\,Wille.  \newblock \textsl{Boolean Concept Logic} in B. Ganter \& G.W. Mineau (Eds.) ICCS 2000 \newblock \textsl{Conceptual Structures: Logical, Linguistic, and Computational Issues, 8th International Conference on Conceptual Structures, ICCS 2000, Darmstadt, Germany, August 14-18, 2000, Proceedings}. LNAI~\textbf{1867} Springer (2000) 317-331.
\bibitem[Wi04]{Wi04} R.\,Wille. \newblock \textsl{Preconcept Algebras and Generalized Double Boolean Algebras}. In P.\,Eklund (Ed.) \textsl{Concept Lattices, Second International Conference on Formal Concept Analysis, ICFCA 2004, Sydney, Australia, February 23-26, 2004, Proceedings}. LNAI~\textbf{2961} Springer (2004) 1-13.
\end{thebibliography}
\end{document}